\documentclass[12pt,a4paper]{article}
\pdfoutput=1 

\usepackage[margin=1in,includefoot]{geometry}
\usepackage[
  bookmarks=true,
  bookmarksnumbered=true,
  bookmarksopen=true,
  pdfborder={0 0 0},
  breaklinks=true,
  colorlinks=true,
  linkcolor=black,
  citecolor=black,
  filecolor=black,
  urlcolor=black,
]{hyperref}

\usepackage{amssymb,amsmath,amsthm}
\usepackage[capitalise]{cleveref}

\usepackage{enumitem}
\usepackage{subfigure}
\usepackage{pgfplots}
\usepackage{mathtools} 

\usepackage{fnbreak} 
\allowdisplaybreaks  

\theoremstyle{plain}
\newtheorem{lemma}{Proposition}
\crefname{lemma}{Proposition}{Propositions}
\newtheorem{cor}[lemma]{Corollary}
\newtheorem{remark}[lemma]{Remark}

\theoremstyle{definition}
\newtheorem{defn}[lemma]{Definition}
\newtheorem{ex}[lemma]{Example}

\newlist{parts}{enumerate}{1}
\crefname{partsi}{Part}{Parts}
\setlist[parts,1]{label=\roman*.,ref=\roman*}

\newcommand{\eqdef}{\triangleq}

\newcommand{\eqby}[1]{\underset{\mbox{\scriptsize #1}}{\scalebox{2.5}[1]{=}}}

\newcommand{\prob}[1]{\mathbf{Pr}\bigl[#1\bigr]}
\newcommand{\probeq}[2]{\prob{#1\!=\!#2}}
\newcommand{\expect}[1]{\mathbf{E}\bigl[#1\bigr]}
\newcommand{\expectsqr}[1]{\mathbf{E}^2\bigl[#1\bigr]}
\newcommand{\variance}[1]{\mathbf{Var}\bigl[#1\bigr]}

\newcommand{\probspace}{\bigl(S_n,2^{S_n},U(S_n)\bigr)}
\newcommand{\spancyc}[1]{\spancycrv{#1}(\pi)}
\newcommand{\spanlen}[1]{\spanlenrv{#1}(\pi)}
\newcommand{\spancycrv}[1]{C^n_{#1}}
\newcommand{\spanlenrv}[1]{L^n_{#1}}

\newcommand{\probnml}{\probeq{\spanlenrv{m}}{\ell}}
\newcommand{\pnml}[3]{p^{#1}_{#2}(#3)}
\newcommand{\nnml}[3]{\Pi^{#1}_{#2}(#3)}

\newcommand{\lfollowing}[1]{\ell^{\pi,#1}_j}

\hypersetup{
  pdfauthor	= {Yannai A. Gonczarowski <yannai@gonch.name>},
  pdftitle	= {"Secret Santa" and the Combined Length of Spanned Cycles in a Random Permutation},
}

\title{``Secret Santa'' and the Combined Length of \\ Spanned Cycles in a Random Permutation}
\author{Yannai A. Gonczarowski\thanks{
Einstein Institute of Mathematics, Rachel \& Selim Benin School of Computer Science and Engineering, and Federmann Center for the Study of Rationality, The Hebrew University of Jerusalem, Israel; and Microsoft Research. \emph{Email}: \href{mailto:yannai@gonch.name}{yannai@gonch.name}.}}
\date{December 25, 2014}

\begin{document}

\maketitle

\begin{abstract}
In many schools throughout the United States, it is customary to hold a yearly ``Secret Santa'' event. With the approach of Christmas, all the students place their names in a hat; each student, in turn, draws a name out of the hat, thus becoming \emph{Secret Santa} to the student whose name is drawn. On
the day before Christmas break, each student places a gift in front of the locker of the student to whom she or he is Secret Santa.

Unfortunately, children being children, some students might forget to bring gifts on this day.
Students who do not find gifts next to their lockers, feeling cheated, angrily take their gifts back from their recipients, who in turn reclaim their gifts from \emph{their} recipients, and so forth.

We analyze the distribution of the number $L$ of students whose Christmas is ruined, as a function of the number of students who do not bring gifts.
We give a simple, explicit formula for the probability of every possible value for $L$ (backed by three proofs of distinct flavors),
as well as closed-form formulae for its expectation and variance.
Notably, we show that if $m$ kids forget to bring gifts, then the expected fraction of kids whose Christmas is \emph{not} ruined is less than $\frac{1}{m+1}$ (regardless of the total number of students), with low probability for a large deviation from this fraction.

The underlying theoretical results are applicable to the study of manipulation in matching markets within game theory.
\end{abstract}

\paragraph{Setting.}

In many schools throughout the United States, it is customary to hold a yearly ``Secret Santa'' event.\footnote{Similar traditions exist in other countries; see, e.g., \url{http://en.wikipedia.org/wiki/Secret\_Santa}.} With the approach of Christmas, all the students place their names in a hat; each student, in turn, draws a name out of the hat, thus becoming \emph{Secret~Santa} to the student whose name is drawn. On the morning of the day before Christmas break, each student places a gift in front of the locker of the student to whom she or he is Secret Santa. At noon, students approach their lockers and find their gifts.
Let us now assume that, unfortunately, some students forget to bring gifts on this day. Students who do not find gifts next to their lockers, feeling cheated, run to where they left their gifts, and angrily take them back from their recipients, who in turn reclaim their gifts from \emph{their} recipients, and so forth.
We analyze the distribution of the number of students whose Christmas is ruined.

\paragraph{Abstraction.}

As the assignment of Secret Santas to students is a one-to-one mapping from the set of students onto itself, it is a \emph{permutation} of the set of students; a student's Christmas is ruined if and only if some student belonging to the same \emph{cycle} of this permutation does not bring a gift (equivalently, a student's Christmas is not ruined if and only if its cycle is disjoint from the set of forgetful students). More abstractly, therefore, given a random permutation of a fixed finite set of objects (e.g., the set of all students), we are interested in the \emph{combined} length $L$ of \emph{all} cycles of the permutation that intersect a given subset $M$ of these objects (e.g., the students who do not bring gifts) --- as explained, in the above-described setting this is the number of students whose Christmas is ruined.\footnote{We assume for simplicity that students may draw their own name out of the hat, in which case these students are their own Secret Santa and as long as they bring a gift, we regard their Christmas as not ruined. As in expectation there will be only one such student (see, e.g., \cite[p.\ 13]{arratia-permutations}),
the number of such students is negligible compared to the number of students who either do not receive gifts or receive gifts that are subsequently taken away,
which, as we will show, is of the order of magnitude of the total number of students even if only one student forgets to bring a gift.
}
When $M$ consists of a single element (i.e., when only one student does not bring a gift), $L$ is simply the well-studied
length of the cycle that contains that element (for an analysis of this special case see, e.g., \cite[p.\ 24]{arratia-permutations}).
The question of the distribution of $L$ arises naturally also during analysis of the limits of manipulation in matching markets within game theory; for more information,
the interested reader is referred to~\cite{blacklists} (matching markets were first defined in \cite{Gale-Shapley}).

\paragraph{Notation.}

We commence by formally defining the problem at hand.

\begin{defn}
Throughout this paper, we use the following standard notation.
\begin{itemize}
\item
$\mathbb{P}\eqdef\{1,2,3,\ldots\}$ --- the positive integers~\cite{stanley-enumerative-combinatorics}; throughout this paper,
the symbols $k,\ell,\tilde{\ell},m,\tilde{m},n,\tilde{n}$ denote elements of~$\mathbb{P}$.
\item
$[n]\eqdef\{1,2,\ldots,n\}$ --- the positive integers up to  $n$~\cite{stanley-enumerative-combinatorics}.
\item
$[m,n]\eqdef\{m,m+1,\ldots,n\}$ --- the integers from $m$ up to $n$~\cite{stanley-enumerative-combinatorics}. ($[m,n]=\emptyset$ if $m>n$.)
\item
$S_N \eqdef \{ \pi : N \mapsto N \mid \mbox{$\pi$ is a bijection} \}$ --- the set of permutations of a set $N$.
\item
$S_n \eqdef S_{[n]}$ --- the set of permutations of $[n]$.
\item
Furthermore, we denote $n^{(k)}\eqdef n\cdot(n+1)\cdot\ldots\cdot(n+k-1)=\frac{(n+k-1)!}{(n-1)!}$ --- the $k$th rising factorial of $n$.
\end{itemize}
\end{defn}

\begin{defn}[Spanned Cycles]
Let $n \in \mathbb{P}$ and $\pi \in S_n$.
For every $M \subseteq [n]$,
we define
\[
\spancyc{M} \:\eqdef\: \bigcup_{m \in M}\bigl\{\pi^{\ell}(m) \mid \ell \in \mathbb{P}
\bigr\} \:\supseteq\: M,
\]
the set of all elements of all cycles of $\pi$ that contain at least one element of $M$.
\end{defn}

Given $n$ and $M$, we study the distribution of $\bigl|\spancyc{M}\bigr|$, i.e., the combined length of all cycles of $\pi$ that intersect $M$, for a random permutation $\pi$ that is uniformly distributed in~$S_n$.
More formally, in the probability space $\probspace$, consisting of $S_n$ as sample space and with the uniform measure over possible outcomes,
we study the distribution of the random variable $\bigl|\spancycrv{M}\bigr|$;
henceforth we work in this space, and denote the outcome of the experiment underlying it by $\pi \in S_n$.
We note that since $\pi\sim U(S_n)$, the distribution of $\bigl|\spancycrv{M}\bigr|$ is the same for sets
$M\subseteq [n]$ of equal size, i.e., this distribution depends on $M$ only through $|M|$; for ease of presentation, we thus consider
only subsets~$M\subseteq[n]$ of the form $M=[m]$ for some $m \le n$, and
define

\begin{defn}[Combined Spanned-Cycles Length]
$\spanlenrv{m}(\pi)\eqdef\bigl|\spancycrv{[m]}(\pi)\bigr| \in [m,n]$, the combined length of all cycles of $\pi$ that contain at least one element less than or equal to~$m$.
\end{defn}

\paragraph{Results.}

We now state the main result of this paper.

\pagebreak[2]

\begin{lemma}[Distribution of {$\spanlenrv{m}$}]\label{distrib}
Let
$m \le n$.
\begin{parts}
\item\label{prob}
$\probnml = \frac{\binom{\ell-1}{m-1}}{\binom{n}{m}}$, for
all $\ell \in [m,n]$.
\item\label{expect}
$\expect{\spanlenrv{m}}=\frac{m \cdot (n+1)}{m+1}$.
\item\label{moments}
$\expect{{\spanlenrv{m}}^{(k)}}=\frac{m\cdot\left[(n+1)^{(k)}\right]}{m+k}$, for
all $k \in \mathbb{P}$.
\item\label{variance}
$\variance{\spanlenrv{m}}=\frac{m\cdot(n+1)\cdot(n-m)}{(m+1)^2\cdot(m+2)}$.
\end{parts}
\end{lemma}

\begin{remark}[Equivalent Formulations of \cref{distrib}(\labelcref{prob})]
\leavevmode
\begin{itemize}
\item
$\probnml = \frac{m}{n}\cdot\prod_{j=1}^{m-1}\frac{\ell-j}{n-j}$.
\item
$\bigl(\probnml\bigr)_{\ell=m}^n\!$ is the prefix of length $n\!-\!m\!+\!1\!$ of
the $m$th diagonal\footnote{The sequences known nowadays
as diagonals of Pascal's triangle are depicted as rows and columns in Pascal's treatise.}\!
of Pascal's triangle~\cite{pascal-triangle},
normalized to sum-up to~1.
\end{itemize}
\end{remark}

\begin{cor}\label{frac}
The expected fraction of the elements of~$[n]$ that are contained in cycles of~$\pi$ that are disjoint from $[m]$ is less than
$\frac{1}{m+1}$, \emph{regardless of the value of $n$}. Furthermore, the standard deviation of this fraction is less than $\frac{1}{m+1}$ as well.
\end{cor}

\begin{cor}\label{full}
$\probeq{\spancycrv{[m]}}{[n]}=\frac{m}{n}$.
\end{cor}

\cref{frac} shows that as $m$ grows, $\spancycrv{[m]}$ quickly grows, \emph{regardless of $n$}, to cover almost all of~$[n]$,
and its size $\spanlenrv{m}$ concentrates on large values (see also \cref{fig}); nonetheless, \cref{full} shows that the probability for $\spancycrv{[m]}$ to cover all of $[n]$ grows considerably slower in a sense, esp.\ for large $n$.
This is demonstrated by the following example.

\newcommand{\plotnml}[3]{
\subfigure[$m=#1$]{
\begin{tikzpicture}
\begin{axis}[
    bar width=1,
    ybar,ymin=0,ymax=1,
    width=8.12cm,
    height=3.82cm,
    tickpos=left,ymajorgrids=true]
\addplot
coordinates
{#2};
\addplot[red,sharp plot,update limits=false,style=dashed]
coordinates {(#3,0) (#3,1)};
\end{axis}
\end{tikzpicture}
}
}

\begin{figure}
\begin{center}
\plotnml{1}{(1,0.01) (2,0.01) (3,0.01) (4,0.01) (5,0.01) (6,0.01) (7,0.01) (8,0.01) (9,0.01) (10,0.01) (11,0.01) (12,0.01) (13,0.01) (14,0.01) (15,0.01) (16,0.01) (17,0.01) (18,0.01) (19,0.01) (20,0.01) (21,0.01) (22,0.01) (23,0.01) (24,0.01) (25,0.01) (26,0.01) (27,0.01) (28,0.01) (29,0.01) (30,0.01) (31,0.01) (32,0.01) (33,0.01) (34,0.01) (35,0.01) (36,0.01) (37,0.01) (38,0.01) (39,0.01) (40,0.01) (41,0.01) (42,0.01) (43,0.01) (44,0.01) (45,0.01) (46,0.01) (47,0.01) (48,0.01) (49,0.01) (50,0.01) (51,0.01) (52,0.01) (53,0.01) (54,0.01) (55,0.01) (56,0.01) (57,0.01) (58,0.01) (59,0.01) (60,0.01) (61,0.01) (62,0.01) (63,0.01) (64,0.01) (65,0.01) (66,0.01) (67,0.01) (68,0.01) (69,0.01) (70,0.01) (71,0.01) (72,0.01) (73,0.01) (74,0.01) (75,0.01) (76,0.01) (77,0.01) (78,0.01) (79,0.01) (80,0.01) (81,0.01) (82,0.01) (83,0.01) (84,0.01) (85,0.01) (86,0.01) (87,0.01) (88,0.01) (89,0.01) (90,0.01) (91,0.01) (92,0.01) (93,0.01) (94,0.01) (95,0.01) (96,0.01) (97,0.01) (98,0.01) (99,0.01) (100,0.01)}{50.5}
~
\plotnml{10}{(1,0.0) (2,0.0) (3,0.0) (4,0.0) (5,0.0) (6,0.0) (7,0.0) (8,0.0) (9,0.0) (10,5.776904234533874e-14) (11,5.776904234533874e-13) (12,3.1772973289936307e-12) (13,1.2709189315974523e-11) (14,4.13048652769172e-11) (15,1.1565362277536816e-10) (16,2.891340569384204e-10) (17,6.608778444306752e-10) (18,1.4043654194151848e-9) (19,2.8087308388303696e-9) (20,5.336588593777702e-9) (21,9.702888352323095e-9) (22,1.6980054616565415e-8) (23,2.8735477043418395e-8) (24,4.7208283714187366e-8) (25,7.553325394269978e-8) (26,1.1802070928546842e-7) (27,1.8050226126012816e-7) (28,2.707533918901922e-7) (29,3.9900499857502016e-7) (30,5.785572479337793e-7) (31,8.265103541911132e-7) (32,1.164628226360205e-6) (33,1.6203523149359373e-6) (34,2.2279844330369135e-6) (35,3.030058828930203e-6) (36,4.0789253466368115e-6) (37,5.438567128849082e-6) (38,7.186677991693429e-6) (39,9.417026333943115e-6) (40,1.2242134234126048e-5) (41,1.5796302237582e-5) (42,2.0239012241901937e-5) (43,2.5758742853329735e-5) (44,3.25772336086229e-5) (45,4.0954236536554505e-5) (46,5.119279567069313e-5) (47,6.364509732032119e-5) (48,7.871893615934463e-5) (49,9.688484450380878e-5) (50,0.00011868393451716575) (51,0.00014473650550873873) (52,0.00017575147097489703) (53,0.00021253666257429405) (54,0.00025601007082812693) (55,0.00030721208499375234) (56,0.00036731879727513865) (57,0.0004376564393065482) (58,0.000519717021676526) (59,0.0006151752501477246) (60,0.0007259067951743151) (61,0.0008540079943227236) (62,0.001001817070263195) (63,0.0011719369501192092) (64,0.0013672597751390775) (65,0.0015909931928891083) (66,0.0018466885274605722) (67,0.002138270926533294) (68,0.002470071587547081) (69,0.002846862168698331) (70,0.0032738914940030803) (71,0.0037569246652494363) (72,0.004302284697301774) (73,0.004916896796916313) (74,0.00560833540898267) (75,0.0063848741579187315) (76,0.00725553881581674) (77,0.00823016343286675) (78,0.009319449769569702) (79,0.010535030174296186) (80,0.011889534053848553) (81,0.013396658088843438) (82,0.015071240349948868) (83,0.016929338475285032) (84,0.018988312073630508) (85,0.021266909522466168) (86,0.02378535933433716) (87,0.026565466269519427) (88,0.0296307123775409) (89,0.0330063631547291) (90,0.03671957900963613) (91,0.04079953223292903) (92,0.04527752967312856) (93,0.05018714132443165) (94,0.05556433503776362) (95,0.06144761757117388) (96,0.06787818220071534) (97,0.07490006311803071) (98,0.08256029684601113) (99,0.09090909090909091) (100,0.1)}{91.81818181818181}

\plotnml{20}{(1,0.0) (2,0.0) (3,0.0) (4,0.0) (5,0.0) (6,0.0) (7,0.0) (8,0.0) (9,0.0) (10,0.0) (11,0.0) (12,0.0) (13,0.0) (14,0.0) (15,0.0) (16,0.0) (17,0.0) (18,0.0) (19,0.0) (20,1.8657295267325187e-21) (21,3.731459053465037e-20) (22,3.918032006138288e-19) (23,2.8732234711680787e-18) (24,1.6521034959216452e-17) (25,7.930096780423896e-17) (26,3.30420699184329e-16) (27,1.2272768826846507e-15) (28,4.142059479060696e-15) (29,1.2886407268188833e-14) (30,3.737058107774761e-14) (31,1.0191976657567531e-13) (32,2.6329273032049455e-13) (33,6.481051823273711e-13) (34,1.527676501200232e-12) (35,3.462733402720526e-12) (36,7.574729318451151e-12) (37,1.6040603262602437e-11) (38,3.297235115090501e-11) (39,6.594470230181002e-11) (40,1.2859216948852953e-10) (41,2.449374656924372e-10) (42,4.5647436788136024e-10) (43,8.335618891746579e-10) (44,1.4934650514379287e-9) (45,2.6284984905307545e-9) (46,4.549324310533998e-9) (47,7.75070067720607e-9) (48,1.3010104708167334e-8) (49,2.1533966413518343e-8) (50,3.5172145142079965e-8) (51,5.672926635819349e-8) (52,9.041226825837087e-8) (53,1.424678166495541e-7) (54,2.220821847772461e-7) (55,3.42641085084894e-7) (56,5.234794355463658e-7) (57,7.92293199745851e-7) (58,1.1884397996187764e-6) (59,1.7674232917407444e-6) (60,2.606949355317598e-6) (61,3.815047837050143e-6) (62,5.54090281095378e-6) (63,7.989208704165914e-6) (64,1.1439094280964832e-5) (65,1.6268934088483317e-5) (66,2.2988711211987297e-5) (67,3.228201999981195e-5) (68,4.506031958307085e-5) (69,6.253268840099627e-5) (70,8.629510999337485e-5) (71,0.00011844426861835765) (72,0.00016172198215198832) (73,0.00021969778707439923) (74,0.0002969988603042805) (75,0.0003995984665912137) (76,0.0005351765177560898) (77,0.0007135686903414531) (78,0.0009473239509705497) (79,0.0012523943758593707) (80,0.0016489859282148383) (81,0.00216260449601946) (82,0.002825338131896391) (83,0.0036774242351667315) (84,0.004769159554981855) (85,0.006163221578745782) (86,0.007937482336263508) (87,0.010188410162965098) (88,0.013035171826146523) (89,0.016624566966679622) (90,0.021136949429064092) (91,0.026793316177686877) (92,0.03386377461346536) (93,0.0426776337594358) (94,0.0536354045895612) (95,0.0672230404189167) (96,0.08402880052364589) (97,0.10476317987363643) (98,0.1302824159967017) (99,0.16161616161616163) (100,0.2)}{96.19047619047619}
~
\plotnml{30}{(1,0.0) (2,0.0) (3,0.0) (4,0.0) (5,0.0) (6,0.0) (7,0.0) (8,0.0) (9,0.0) (10,0.0) (11,0.0) (12,0.0) (13,0.0) (14,0.0) (15,0.0) (16,0.0) (17,0.0) (18,0.0) (19,0.0) (20,0.0) (21,0.0) (22,0.0) (23,0.0) (24,0.0) (25,0.0) (26,0.0) (27,0.0) (28,0.0) (29,0.0) (30,3.4045636339269157e-26) (31,1.0213690901780745e-24) (32,1.5831220897760155e-23) (33,1.68866356242775e-22) (34,1.3931474390028939e-21) (35,9.473402585219677e-21) (36,5.526151508044813e-20) (37,2.842020775565903e-19) (38,1.3144346086992301e-18) (39,5.54983501450786e-18) (40,2.164435655658066e-17) (41,7.870675111483875e-17) (42,2.689147329756991e-16) (43,8.688014449984124e-16) (44,2.6684615810665523e-15) (45,7.827487304461887e-15) (46,2.2014808043799058e-14) (47,5.956948058910334e-14) (48,1.5554253264932536e-13) (49,3.9294955616671674e-13) (50,9.62726412608456e-13) (51,2.2922057443058473e-12) (52,5.3137496799817375e-12) (53,1.2013694928654363e-11) (54,2.6530242967445052e-11) (55,5.730532480968131e-11) (56,1.2122280248201815e-10) (57,2.5142507181455617e-10) (58,5.118296104796322e-10) (59,1.0236592209592644e-9) (60,2.0131964678865534e-9) (61,3.896509292683652e-9) (62,7.427720839178211e-9) (63,1.3955111879668154e-8) (64,2.5858001424090993e-8) (65,4.72832026040521e-8) (66,8.537244914620519e-8) (67,1.5228599036890655e-7) (68,2.6850424617675624e-7) (69,4.6816124974408785e-7) (70,8.075781558085516e-7) (71,1.3787919733316734e-6) (72,2.3308150025368764e-6) (73,3.902760004247793e-6) (74,6.475033643411111e-6) (75,1.0647833102498271e-5) (76,1.736059744972544e-5) (77,2.8072455450619864e-5) (78,4.503289728536936e-5) (79,7.168502016854716e-5) (80,0.0001132623318663045) (81,0.00017766640292753649) (82,0.00027674958917558567) (83,0.0004281786096678873) (84,0.0006581263815265675) (85,0.001005138473604212) (86,0.0015256566117206791) (87,0.0023018678703154104) (88,0.003452801805473116) (89,0.005149941675959902) (90,0.007639080152673854) (91,0.01127077399574831) (92,0.016542587638920907) (93,0.02415742956794799) (94,0.03510376484092442) (95,0.05076544453918301) (96,0.07307147320033917) (97,0.10469942428705316) (98,0.14935064935064934) (99,0.21212121212121213) (100,0.3)}{97.74193548387096}

\plotnml{40}{(1,0.0) (2,0.0) (3,0.0) (4,0.0) (5,0.0) (6,0.0) (7,0.0) (8,0.0) (9,0.0) (10,0.0) (11,0.0) (12,0.0) (13,0.0) (14,0.0) (15,0.0) (16,0.0) (17,0.0) (18,0.0) (19,0.0) (20,0.0) (21,0.0) (22,0.0) (23,0.0) (24,0.0) (25,0.0) (26,0.0) (27,0.0) (28,0.0) (29,0.0) (30,0.0) (31,0.0) (32,0.0) (33,0.0) (34,0.0) (35,0.0) (36,0.0) (37,0.0) (38,0.0) (39,0.0) (40,7.274719675172519e-29) (41,2.9098878700690075e-27) (42,5.965270133641464e-26) (43,8.35137818709805e-25) (44,8.977731551130404e-24) (45,7.900403764994756e-23) (46,5.925302823746067e-22) (47,3.89377042703313e-21) (48,2.2875901258819637e-20) (49,1.2200480671370472e-19) (50,5.978235528971532e-19) (51,2.7173797858961507e-18) (52,1.154886409005864e-17) (53,4.619545636023456e-17) (54,1.7488279907803085e-16) (55,6.29578076680911e-16) (56,2.1641746385906315e-15) (57,7.129045868298552e-15) (58,2.2575311916278746e-14) (59,6.891411006021933e-14) (60,2.0329662467764702e-13) (61,5.808474990789915e-13) (62,1.6105317019917491e-12) (63,4.341433283629932e-12) (64,1.1396262369528573e-11) (65,2.917443166599315e-11) (66,7.293607916498287e-11) (67,1.7828819351440255e-10) (68,4.266181773380347e-10) (69,1.0003460709995297e-9) (70,2.3007959632989182e-9) (71,5.195345723578202e-9) (72,1.1527173324189137e-8) (73,2.5150196343685388e-8) (74,5.3998950973206866e-8) (75,1.1416921062906593e-7) (76,2.3785252214388738e-7) (77,4.885619373766335e-7) (78,9.89980767842126e-7) (79,1.979961535684252e-6) (80,3.910424032976397e-6) (81,7.630095674100286e-6) (82,1.4715184514336268e-5) (83,2.8061514655245905e-5) (84,5.293422082694114e-5) (85,9.881054554362346e-5) (86,0.00018258470372191292) (87,0.0003340911600017981) (88,0.0006055402275032591) (89,0.0010875008167405469) (90,0.0019357514537981734) (91,0.0034160319772908944) (92,0.005978055960259065) (93,0.010377002798940264) (94,0.01787150482039712) (95,0.03054402642031508) (96,0.05181575910589165) (97,0.08726864691518595) (98,0.14594928880643165) (99,0.24242424242424243) (100,0.4)}{98.53658536585365}
~
\plotnml{50}{(1,0.0) (2,0.0) (3,0.0) (4,0.0) (5,0.0) (6,0.0) (7,0.0) (8,0.0) (9,0.0) (10,0.0) (11,0.0) (12,0.0) (13,0.0) (14,0.0) (15,0.0) (16,0.0) (17,0.0) (18,0.0) (19,0.0) (20,0.0) (21,0.0) (22,0.0) (23,0.0) (24,0.0) (25,0.0) (26,0.0) (27,0.0) (28,0.0) (29,0.0) (30,0.0) (31,0.0) (32,0.0) (33,0.0) (34,0.0) (35,0.0) (36,0.0) (37,0.0) (38,0.0) (39,0.0) (40,0.0) (41,0.0) (42,0.0) (43,0.0) (44,0.0) (45,0.0) (46,0.0) (47,0.0) (48,0.0) (49,0.0) (50,9.911653021418339e-30) (51,4.955826510709169e-28) (52,1.2637357602308382e-26) (53,2.190475317733453e-25) (54,2.902379795996825e-24) (55,3.1345701796765706e-23) (56,2.8733559980368566e-22) (57,2.2986847984294852e-21) (58,1.6378129188810086e-20) (59,1.0554794366122054e-19) (60,6.227328676012013e-19) (61,3.3967247323701883e-18) (62,1.7266684056215124e-17) (63,8.234880088348752e-17) (64,3.7056960397569384e-16) (65,1.5810969769629602e-15) (66,6.423206468912026e-15) (67,2.4937154526364337e-14) (68,9.282163073702281e-14) (69,3.3220373105881846e-13) (70,1.1461028721529237e-12) (71,3.8203429071764125e-12) (72,1.232928847316024e-11) (73,3.8596033481197276e-11) (74,1.1739626850530837e-10) (75,3.474929547757128e-10) (76,1.0023835233914792e-9) (77,2.8215239917686083e-9) (78,7.759190977363673e-9) (79,2.0869548146012637e-8) (80,5.4956476784499945e-8) (81,1.418231658954837e-7) (82,3.5898988867294314e-7) (83,8.920354809448891e-7) (84,2.1776160270125233e-6) (85,5.226278464830056e-6) (86,1.2339824153070965e-5) (87,2.868175343686765e-5) (88,6.56661197107233e-5) (89,0.00014816970601393977) (90,0.000329677595881016) (91,0.0007236825275436936) (92,0.0015679788096780028) (93,0.0033547453602413083) (94,0.007090711784146402) (95,0.014811709060216928) (96,0.030589399146100177) (97,0.062480474851608875) (98,0.12626262626262627) (99,0.25252525252525254) (100,0.5)}{99.01960784313725}

\plotnml{60}{(1,0.0) (2,0.0) (3,0.0) (4,0.0) (5,0.0) (6,0.0) (7,0.0) (8,0.0) (9,0.0) (10,0.0) (11,0.0) (12,0.0) (13,0.0) (14,0.0) (15,0.0) (16,0.0) (17,0.0) (18,0.0) (19,0.0) (20,0.0) (21,0.0) (22,0.0) (23,0.0) (24,0.0) (25,0.0) (26,0.0) (27,0.0) (28,0.0) (29,0.0) (30,0.0) (31,0.0) (32,0.0) (33,0.0) (34,0.0) (35,0.0) (36,0.0) (37,0.0) (38,0.0) (39,0.0) (40,0.0) (41,0.0) (42,0.0) (43,0.0) (44,0.0) (45,0.0) (46,0.0) (47,0.0) (48,0.0) (49,0.0) (50,0.0) (51,0.0) (52,0.0) (53,0.0) (54,0.0) (55,0.0) (56,0.0) (57,0.0) (58,0.0) (59,0.0) (60,7.274719675172519e-29) (61,4.364831805103511e-27) (62,1.3312737005565707e-25) (63,2.7512989811502462e-24) (64,4.333295895311638e-23) (65,5.546618745998896e-22) (66,6.0088369748321374e-21) (67,5.665474861984587e-20) (68,4.744835196912092e-19) (69,3.584986593222469e-18) (70,2.4736407493235036e-17) (71,1.574135022296775e-16) (72,9.31363221525592e-16) (73,5.158319380757125e-15) (74,2.6896951056805007e-14) (75,1.3269162521357135e-13) (76,6.219919931886158e-13) (77,2.7806700871961648e-12) (78,1.189508870633915e-11) (79,4.8832469426023876e-11) (80,1.928882542327943e-10) (81,7.348123970773116e-10) (82,2.7054456437846474e-9) (83,9.645501860449612e-9) (84,3.335736060072158e-8) (85,1.120807316184245e-7) (86,3.664177764448493e-7) (87,1.1671084731206312e-6) (88,3.6263727557676754e-6) (89,1.1004165603708808e-5) (90,3.26456912910028e-5) (91,9.477781342549199e-5) (92,0.00026952440692874286) (93,0.0007514013768922528) (94,0.00205530376620528) (95,0.005519958686379894) (96,0.01456655764461361) (97,0.03779431172656504) (98,0.09647495361781076) (99,0.24242424242424243) (100,0.6)}{99.34426229508196}
~
\plotnml{70}{(1,0.0) (2,0.0) (3,0.0) (4,0.0) (5,0.0) (6,0.0) (7,0.0) (8,0.0) (9,0.0) (10,0.0) (11,0.0) (12,0.0) (13,0.0) (14,0.0) (15,0.0) (16,0.0) (17,0.0) (18,0.0) (19,0.0) (20,0.0) (21,0.0) (22,0.0) (23,0.0) (24,0.0) (25,0.0) (26,0.0) (27,0.0) (28,0.0) (29,0.0) (30,0.0) (31,0.0) (32,0.0) (33,0.0) (34,0.0) (35,0.0) (36,0.0) (37,0.0) (38,0.0) (39,0.0) (40,0.0) (41,0.0) (42,0.0) (43,0.0) (44,0.0) (45,0.0) (46,0.0) (47,0.0) (48,0.0) (49,0.0) (50,0.0) (51,0.0) (52,0.0) (53,0.0) (54,0.0) (55,0.0) (56,0.0) (57,0.0) (58,0.0) (59,0.0) (60,0.0) (61,0.0) (62,0.0) (63,0.0) (64,0.0) (65,0.0) (66,0.0) (67,0.0) (68,0.0) (69,0.0) (70,3.4045636339269157e-26) (71,2.3831945437488406e-24) (72,8.460340630308385e-23) (73,2.030481751274012e-21) (74,3.7056291960750725e-20) (75,5.484331210191107e-19) (76,6.855414012738884e-18) (77,7.443020928116503e-17) (78,7.163907643312133e-16) (79,6.208719957537182e-15) (80,4.904888766454374e-14) (81,3.567191830148636e-13) (82,2.407854485350329e-12) (83,1.5188005215286692e-11) (84,9.004317377634253e-11) (85,5.042417731475181e-10) (86,2.6787844198461903e-9) (87,1.355149765333955e-8) (88,6.54989053244745e-8) (89,3.0336335097651343e-7) (90,1.349966911845485e-6) (91,5.7855724793377924e-6) (92,2.3931231619079048e-5) (93,9.572492647631619e-5) (94,0.00037093409009572527) (95,0.001394712178759927) (96,0.0050960637300843485) (97,0.01811933770696657) (98,0.06277056277056277) (99,0.21212121212121213) (100,0.7)}{99.5774647887324}

\plotnml{80}{(1,0.0) (2,0.0) (3,0.0) (4,0.0) (5,0.0) (6,0.0) (7,0.0) (8,0.0) (9,0.0) (10,0.0) (11,0.0) (12,0.0) (13,0.0) (14,0.0) (15,0.0) (16,0.0) (17,0.0) (18,0.0) (19,0.0) (20,0.0) (21,0.0) (22,0.0) (23,0.0) (24,0.0) (25,0.0) (26,0.0) (27,0.0) (28,0.0) (29,0.0) (30,0.0) (31,0.0) (32,0.0) (33,0.0) (34,0.0) (35,0.0) (36,0.0) (37,0.0) (38,0.0) (39,0.0) (40,0.0) (41,0.0) (42,0.0) (43,0.0) (44,0.0) (45,0.0) (46,0.0) (47,0.0) (48,0.0) (49,0.0) (50,0.0) (51,0.0) (52,0.0) (53,0.0) (54,0.0) (55,0.0) (56,0.0) (57,0.0) (58,0.0) (59,0.0) (60,0.0) (61,0.0) (62,0.0) (63,0.0) (64,0.0) (65,0.0) (66,0.0) (67,0.0) (68,0.0) (69,0.0) (70,0.0) (71,0.0) (72,0.0) (73,0.0) (74,0.0) (75,0.0) (76,0.0) (77,0.0) (78,0.0) (79,0.0) (80,1.8657295267325187e-21) (81,1.4925836213860148e-19) (82,6.04496366661336e-18) (83,1.6522900688743184e-16) (84,3.4285018929142108e-15) (85,5.759883180095874e-14) (86,8.159834505135821e-13) (87,1.0024939534881153e-11) (88,1.0902121744183253e-10) (89,1.0659852372090292e-9) (90,9.48726861116036e-9) (91,7.762310681858476e-8) (92,5.886418933742677e-7) (93,4.165773399264049e-6) (94,2.7672637580825465e-5) (95,0.00017341519550650625) (96,0.0010296527233198808) (97,0.005814509496394621) (98,0.0313337456194599) (99,0.16161616161616163) (100,0.8)}{99.75308641975309}
~
\plotnml{90}{(1,0.0) (2,0.0) (3,0.0) (4,0.0) (5,0.0) (6,0.0) (7,0.0) (8,0.0) (9,0.0) (10,0.0) (11,0.0) (12,0.0) (13,0.0) (14,0.0) (15,0.0) (16,0.0) (17,0.0) (18,0.0) (19,0.0) (20,0.0) (21,0.0) (22,0.0) (23,0.0) (24,0.0) (25,0.0) (26,0.0) (27,0.0) (28,0.0) (29,0.0) (30,0.0) (31,0.0) (32,0.0) (33,0.0) (34,0.0) (35,0.0) (36,0.0) (37,0.0) (38,0.0) (39,0.0) (40,0.0) (41,0.0) (42,0.0) (43,0.0) (44,0.0) (45,0.0) (46,0.0) (47,0.0) (48,0.0) (49,0.0) (50,0.0) (51,0.0) (52,0.0) (53,0.0) (54,0.0) (55,0.0) (56,0.0) (57,0.0) (58,0.0) (59,0.0) (60,0.0) (61,0.0) (62,0.0) (63,0.0) (64,0.0) (65,0.0) (66,0.0) (67,0.0) (68,0.0) (69,0.0) (70,0.0) (71,0.0) (72,0.0) (73,0.0) (74,0.0) (75,0.0) (76,0.0) (77,0.0) (78,0.0) (79,0.0) (80,0.0) (81,0.0) (82,0.0) (83,0.0) (84,0.0) (85,0.0) (86,0.0) (87,0.0) (88,0.0) (89,0.0) (90,5.776904234533874e-14) (91,5.199213811080486e-12) (92,2.3656422840416215e-10) (93,7.254636337727639e-9) (94,1.686702948521676e-7) (95,3.171001543220751e-6) (96,5.020752443432856e-5) (97,0.000688560335099363) (98,0.008348794063079777) (99,0.09090909090909091) (100,0.9)}{99.89010989010988}

\caption{Distribution of $\spanlenrv{m}$ for $n=100$ and varying values of $m$. The (red) dashed vertical line marks the expectation.
As can be seen, as $m$ grows, $\spanlenrv{m}$ quickly concentrates on large values close to $n$.}
\label{fig}
\end{center}
\end{figure}

\begin{ex}
For $n=1000$, $C^{1000}_{[100]}$ is expected to cover more than
$99\%$ of all elements (with standard deviation\ $\sigma<1\%$), while the probability that it covers all $1000$ elements is $\frac{1}{10}$.
\end{ex}

\paragraph{Proofs.}

We present three proofs for \cref{distrib}(\labelcref{prob}). The first proof is recursive, calculating the distribution for $m$ given the distribution for $m\!-\!1$. The second proof is enumerative, directly and succinctly proving the special case in which $\ell=n$ (i.e., \cref{full}), and proving the general case by reduction to this special case. The third proof, also enumerative, provides an interpretation of the nominator and the denominator of the r.h.s.\ of the equality in \cref{distrib}(\labelcref{prob}).

\begin{proof}[\textbf{Probabilistic proof of \cref{distrib}(\labelcref{prob})}]
For all
$m \le \ell \le n$,
we define
\[\pnml{n}{m}{\ell}\eqdef\probnml.\]

Throughout this proof, we make extensive use of the following well-known (see, e.g., \cite[p.\ 7]{riordan-identities}) identity:
\begin{equation}\label{sum}
\forall m,n \in \mathbb{P}:
\binom{n}{m}=\smashoperator[r]{\sum_{j=m-1}^{n-1}}\ \binom{j}{m-1},
\end{equation}
obtained either inductively as in \cite{riordan-identities}, or by conditioning
upon the maximum element in the chosen set of $m$-out-of-$n$ elements.

We prove, by induction on $\tilde{m}$, that
$\pnml{\tilde{n}}{\tilde{m}}{\tilde{\ell}}=\frac{\binom{\tilde{\ell}-1}{\tilde{m}-1}}{\binom{\tilde{n}}{\tilde{m}}}$
for all
$\tilde{m} \le \tilde{\ell} \le \tilde{n}$.

Base: Let $\ell \le n$. We show that the claim
holds for $\tilde{m}=1, \tilde{\ell}=\ell$, and $\tilde{n}=n$.
We observe that $\pnml{n}{1}{\ell}$ is simply the probability that the cycle
of $\pi\sim U(S_n)$
containing the element
$1$ has length $\ell$. It is well established~\cite[p.\ 24]{arratia-permutations}
that the length of this cycle is uniformly distributed in $[n]$, yielding
$\pnml{n}{1}{\ell}=\frac{1}{n}=\frac{\binom{\ell-1}{0}}{\binom{n}{1}}$, as required.

Step: Let
$1<m\le \ell\le n$, and assume that the claim holds
for $\tilde{n}=n$, $\tilde{m}=m\!-\!1$, and all $\tilde{\ell} \in [m-1,n]$;
furthermore, assume that the base case holds whenever
$\tilde{\ell}\le \tilde{n} < n$.
We claim that the following recurrence relation holds:
\begin{equation}\label{recur}
\pnml{n}{m}{\ell} = \pnml{n}{m-1}{\ell}\cdot\tfrac{\ell-m+1}{n-m+1} +
\smashoperator{\sum_{j=m-1}^{\ell-1}}\ \pnml{n}{m-1}{j} \cdot \tfrac{n-j}{n-m+1} \cdot \pnml{n-j}{1}{\ell-j}.
\end{equation}
We justify \cref{recur} using the law of total probability, by conditioning upon the
value of $j\eqdef \spanlenrv{m-1} \in [m-1,n]$.
If~$j>\ell$, then obviously $\spanlenrv{m}\ge j>\ell$ with probability~$1$.
If~$j=\ell$, then $\spanlenrv{m}=\ell$ iff $m\in \spancycrv{[m-1]}$, which
holds with probability $\frac{\left|\spancycrv{[m-1]}\setminus[m-1]\right|}{\left|[m,n]\right|}=
\frac{\ell-m+1}{n-m+1}$. Otherwise, i.e.,
if $m-1\le j<\ell$, then $\spanlenrv{m}=\ell$ iff both $m\notin \spancycrv{[m-1]}$ and $\bigl|\spancycrv{\{m\}}\bigr|=\ell-j$;
the first condition holds with probability
$\frac{\left|[n]\setminus\spancycrv{[m-1]}\right|}{\left|[m,n]\right|}=
\frac{n-j}{n-m+1}$, and the second (conditioned upon the first) --- with probability
$\pnml{n-j}{1}{\ell-j}$, since $\pi|_{[n]\setminus \spancycrv{[m-1]}}$, given $\spancycrv{[m-1]}$,
is uniformly distributed in $S_{[n]\setminus \spancycrv{[m-1]}}\cong S_{n-j}$.

Plugging the induction hypothesis for $\tilde{m}=m\!-\!1$
and the base case for $\tilde{n}=n\!-\!j$ into \cref{recur}, we obtain
\begin{align*}
\pnml{n}{m}{\ell} &= 
\frac{\binom{\ell-1}{m-2}}{\binom{n}{m-1}}
\cdot\frac{\ell-m+1}{n-m+1} +
\smashoperator{\sum_{j=m-1}^{\ell-1}}\ 
\frac{\binom{j-1}{m-2}}{\binom{n}{m-1}}\cdot\frac{n-j}{n-m+1}\cdot\frac{1}{n-j} = \\
&=\frac{1}{\binom{n}{m-1}\cdot(n-m+1)}\cdot
\left(\binom{\ell-1}{m-2}\cdot(\ell-m+1) +
\smashoperator{\sum_{j=m-1}^{\ell-1}}\ \binom{j-1}{m-2}\right) = \\
&=\frac{1}{\binom{n}{m}\cdot m}\cdot
\left(\binom{\ell-1}{m-2}\cdot(\ell-m+1) +
\smashoperator{\sum_{j=m-1}^{\ell-1}}\ \binom{j-1}{m-2}\right) = \\
&=\frac{1}{\binom{n}{m}\cdot m}\cdot
\left((m-1)\cdot\binom{\ell-1}{m-1} +
\smashoperator{\sum_{j=m-1}^{\ell-1}}\ \binom{j-1}{m-2}\right) \eqby{\cref{sum}} \\
&=\frac{1}{\binom{n}{m}\cdot m}\cdot
\left((m-1)\cdot\binom{\ell-1}{m-1} +
\binom{\ell-1}{m-1}\right) = \\
&=\frac{1}{\binom{n}{m}\cdot m}\cdot m \cdot
\binom{\ell-1}{m-1} =
\frac{\binom{\ell-1}{m-1}}{\binom{n}{m}},
\end{align*}
and the proof by induction is complete.
We note that by \cref{sum}, we immediately verify that indeed
$\sum_{\ell=m}^n \pnml{n}{m}{\ell}=1$ for all
$m \le n$.
\end{proof}
 
As mentioned above, before proceeding to the proofs of the remaining parts of \cref{distrib}, we present two additional,
significantly different, proofs of \cref{distrib}(\labelcref{prob}). Both
of these proofs, while of distinct flavors, make use of the following
definition.

\begin{defn}\label{restrict}
Let $n \in \mathbb{P}$.
For every $\pi \in S_n$ and
$k \le n$,
by a slight abuse of notation
we denote by $\pi|_k \in S_k$
the permutation obtained by inspecting the
cycle-structure representation of $\pi$ and removing all elements of $[k+1,n]$
from it. More formally, for every $j \in [k]$, we define
$\pi|_k(j) \eqdef \pi^{\lfollowing{k}}(j)$, where $\lfollowing{k}$ is the smallest positive
integer s.t.\ $\pi^{\lfollowing{k}}(j) \in [k]$.
\end{defn}

\begin{ex}\label{restrict-ex}
If $\pi|_6=(365)(24)(1)$ (in cycle-structure representation), then the
cycle-structure representation of $\pi$ is of the form
\[
\cdots(3\ldots 6\ldots 5\ldots)(2\ldots4\ldots)(1\ldots),
\]
where the first ellipsis stands for zero or more cycles disjoint from the set [6],
and each subsequent
ellipsis stands for zero or more consecutive elements greater than $6$ within a cycle.
(E.g., $\pi\eqdef(8)(3~10~11~6~5~7)(2~4)(1~9) \in S_{11}$ is of this form.)
In fact,
for every $j \in [6]$, the ellipsis immediately following $j$ stands for
precisely $\lfollowing{6}\!-\!1$ (as defined in \cref{restrict}) elements,
while the first ellipsis stands for a product
of cycles of combined length $n\!-\!\sum_{j=1}^6\lfollowing{6}$.
\end{ex}

\begin{proof}[\textbf{Enumerative proof by reduction for \cref{distrib}(\labelcref{prob})}]
For every
$m\le \ell\le n$, we define
\[
\nnml{n}{m}{\ell}\eqdef\bigl\{\pi \in S_n \mid \spanlen{m}\!=\!\ell\bigr\}.
\]
We show that $\bigl|\nnml{n}{m}{\ell}\bigr|=\binom{n-m}{\ell-m}\cdot m\cdot (\ell-1)! \cdot (n-\ell)!=
n!\cdot\frac{\binom{\ell-1}{m-1}}{\binom{n}{m}}$.
We first show this for the special case of $\ell=n$; i.e., we show that
for all
$m \le n$, the set $\nnml{n}{m}{n}$, of permutations
on $[n]$ with all cycles intersecting $[m]$, is of size $m\cdot(n-1)!$.

Consider the following argument for the
equality $|S_n|=n!$, tracing the construction of a permutation $\pi \in S_n$
by iteratively constructing $\pi|_1$, then $\pi|_2$,
and so forth until $\pi|_n=\pi$. Obviously, $\pi|_1=(1)$. To obtain $\pi|_2$ from
$\pi|_1$, a two-way choice is made: the element $2$ may be placed either
(immediately) after $1$ in its cycle, or in a new (singleton) cycle. To obtain $\pi|_3$,
a three-way choice is made: the element $3$ may now
be placed either after $1$ in its cycle, after $2$ in its cycle, or in a new cycle.
More generally, to obtain $\pi|_k$ from $\pi|_{k-1}$, for $k \in [2,n]$,
a $k$-way choice is made: the element~$k$ may be placed either after
some element $j \in [k-1]$ in its cycle (more formally,
setting ${\pi|_k}^{-1}(k)=j$ and $\pi|_k(k)=\pi|_{k-1}(j)$),
or in a new cycle (i.e., having $k$ a fixed point of $\pi|_k$).
Thus, we obtain that there are $n!$ ways to construct a permutation $\pi \in S_n$, each resulting in a distinct outcome (as $\pi$ uniquely determines $\pi|_k$ for all
$k \le n$),
as required. We note that in fact, construction of a permutation $\pi \in \nnml{n}{m}{n}$
may be undertaken in a very similar manner, the only difference being that the elements of $[m+1,n]$
may not be placed in new cycles, thus reducing the choice for each $k \in [m+1,n]$ from a $k$-way choice to a $(k\!-\!1)$-way one. By similar reasoning, we therefore obtain
$\bigl|\nnml{n}{m}{n}\bigr|=m!\cdot m^{(n-m)}=m\cdot(n-1)!$.

We now move on to the general case. A permutation $\pi \in \nnml{n}{m}{\ell}$
may be constructed as follows.
First, choose a subset $I \subseteq [m+1,n]$ of size $\ell\!-\!m$ as the additional
elements, in addition to $[m]$, of $\spancyc{[m]}$. (There are $\binom{n-m}{\ell-m}$ options.)
Next, choose any permutation on $[m] \cup I$ in which all cycles
intersect $[m]$ ---
this permutation constitutes the product of the cycles of $\pi$ that
intersect $[m]$. (There are $m\cdot(\ell-1)!$ options, by the above special case.) Finally, choose any permutation on
$[n] \setminus ([m] \cup I)$ as the product of the
remaining cycles of $\pi$, i.e., those that do not intersect $[m]$. (There are $(n-\ell)!$
options.)
We thus obtain $\bigl|\nnml{n}{m}{\ell}\bigr|=\binom{n-m}{\ell-m}\cdot m\cdot (\ell-1)! \cdot (n-\ell)!$,
as required.
\end{proof}

\sloppy

\begin{proof}[\textbf{Direct enumerative proof for \cref{distrib}(\labelcref{prob})}]
Henceforth, when representing the
cycle structure of any permutation, we write each cycle with its
smallest element first, and write cycles in decreasing order of their
first (i.e., smallest) element.
E.g., the reader may verify that
all cycle-structure representations in \cref{restrict-ex}, and notably
that of the general form (i.e., with ellipses) of $\pi$ in that example,
follow this convention.
It is straightforward to check
(see, e.g., \cite[Section~1.3]{stanley-enumerative-combinatorics}, where
a similar convention is used) that this representation is both unique and
unambiguous even when the parentheses are dispensed with. (Indeed, uniqueness
implies unambiguity, since
the number of ways to order $[n]$ in a row
 equals the number of permutations on~$[n]$.)

\fussy

Let
$m \le n$.
For every $\pi \in S_n$, we denote by
$\pi_{> m}$ the sequence consisting of the elements of $[m+1,n]$,
ordered as in the cycle-structure representation (according to the above
convention) of $\pi$.
We claim that the mapping
$\pi \mapsto \bigl(\pi|_m, (\lfollowing{m})_{j=1}^m, \pi_{> m}\bigr)$ is a bijection
between $S_n$ and $S_m\times \bigl\{(\ell_j)_{j=1}^m\!\in\mathbb{P}^m \mid \sum_{j=1}^m\ell_j\le n\bigr\} \times S_{[m+1,n]}$, where
by a very slight abuse of notation we think of a permutation $\tau \in S_{[m+1,n]}$ as the sequence $\bigl(\tau(m+1),\ldots,\tau(n)\bigr)$.
Under the notation of \cref{restrict-ex}, $\pi|_m$
determines the general form of $\pi$ w.r.t.\ $[m]$, while
$(\lfollowing{m})_{j=1}^m$ determine the number of elements each ellipsis stands
for, and
$\pi_{>m}$, given all of these, determines the exact content of each ellipsis
(the unambiguity of the cycle-structure representation, even when the parentheses are dispensed with,
is used when populating the first ellipsis).
The reader who is not yet convinced of the validity of this bijection
claim,
may verify that this mapping is onto, and that the size of the domain
and the size of the image match (see the last equality of \cref{second-multiplicand} below for the size of the second
multiplicand).

Let $\ell \in [m+1,n]$.
We observe that for every $\pi \in S_n$, by definition
$\spanlen{m}=\sum_{j=1}^m\lfollowing{m}$ (see, e.g., the suffix of \cref{restrict-ex}). Thus, we have that for every
$\sigma \in S_m$ and $\tau \in S_{[m+1,n]}$,
\begin{multline*}
\bigl|\bigl\{\pi \in S_n \;\mid\; \pi|_m\!=\!\sigma \And \pi_{> m}\!=\!\tau
\And \spanlen{m}\!=\!\ell \bigr\}\bigr|=\\*
\textstyle
\bigl|\bigl\{(\ell_j)_{j=1}^m\!\in\mathbb{P}^m \mid \sum_{j=1}^m\ell_j=\ell\bigr\}\bigr|=
\binom{\ell-1}{m-1}.
\end{multline*}
(For the calculation of the number of $m$-compositions of $\ell$ see, e.g., \cite[Section~1.2]{stanley-enumerative-combinatorics}.) For comparison,
dispensing with the conditioning on $\spanlen{m}$ we have
\begin{equation}\label{second-multiplicand}
\textstyle
\bigl|\bigl\{\pi \in S_n \;\mid\; \pi|_m\!=\!\sigma \And \pi_{> m}\!=\!\tau \bigr\}\bigr|=
\textstyle
\bigl|\bigl\{(\ell_j)_{j=1}^m\!\in\mathbb{P}^m \mid \sum_{j=1}^m\ell_j\le n\bigr\}\bigr|=\binom{n}{m},
\end{equation}
since such $(\ell_j)_{j=1}^m$ are in one-to-one correspondence with $(m\!+\!1$)-compositions of ${n\!+\!1}$, where the $(m\!+\!1)$th element designates the successor
of the remainder.
Combining these, we obtain the slightly stronger result that
\[
\probeq{\spanlenrv{m}}{\ell\;\mid\; \pi|_m\!=\!\sigma \And \pi_{> m}\!=\!\tau}=
\frac{\binom{\ell-1}{m-1}}{\binom{n}{m}},
\]
for every choice of $\sigma \in S_m$ and $\tau \in S_{[m+1,n]}$. As the r.h.s.\ depends on neither $\sigma$ nor~$\tau$, we have
\[
\probnml=
\frac{\binom{\ell-1}{m-1}}{\binom{n}{m}},
\]
as required.
\end{proof}

Finally, we prove the remaining parts of \cref{distrib}.

\begin{proof}[\textbf{Proof of \cref{distrib}(\labelcref{expect,moments,variance})}]
We prove \cref{expect}
directly by definition of expectation:
\begin{align*}
\expect{\spanlenrv{m}}
&= \sum_{\ell=m}^n \probnml\cdot \ell =
\frac{1}{\binom{n}{m}}\cdot\left(\sum_{\ell=m}^n \binom{\ell-1}{m-1} \cdot \ell\right) = \\
&= \frac{1}{\binom{n}{m}}\cdot\left(m \cdot \sum_{\ell=m}^n \binom{\ell}{m} \right) \eqby{\cref{sum}} \\
&= \frac{1}{\binom{n}{m}}\cdot\left(m \cdot \binom{n+1}{m+1} \right) =
\frac{m\cdot(n+1)}{m+1}.
\end{align*}
More generally, all \emph{rising-factorial moments} may be calculated in a similar manner:
\begin{align*}
\expect{{\spanlenrv{m}}^{(k)}}
&= \sum_{\ell=m}^n \probnml\cdot \ell^{(k)} =
\frac{1}{\binom{n}{m}}\cdot\left(\sum_{\ell=m}^n \binom{\ell-1}{m-1} \cdot \ell^{(k)}\right) = \\
&= \frac{1}{\binom{n}{m}}\cdot\left(m^{(k)} \cdot \sum_{\ell=m}^n \binom{\ell+k-1}{m+k-1} \right) \eqby{\cref{sum}} \\
&= \frac{1}{\binom{n}{m}}\cdot\left(m^{(k)} \cdot \binom{n+k}{m+k} \right) =
\frac{m\cdot\bigl[(n+1)^{(k)}\bigr]}{m+k}.
\end{align*}
The rising-factorial moments give rise to calculation of the raw moments
and the central moments.
The second raw moment, for instance, is given by
\begin{align*}
\expect{{\spanlenrv{m}}^2}
&= \expect{{\spanlenrv{m}}^{(2)}}-\expect{\spanlenrv{m}}=\\
&= \frac{m\cdot(n+1)\cdot(n+2)}{m+2}-\frac{m\cdot(n+1)}{m+1}=\\
&= m\cdot(n+1)\cdot\left(\frac{n+2}{m+2}-\frac{1}{m+1}\right)=\\
&= m\cdot(n+1)\cdot\frac{(m+1)(n+2)-(m+2)}{(m+2)(m+1)}=\\
&= \frac{m\cdot(n+1)\cdot (mn+n+m)}{(m+2)(m+1)},
\end{align*}
and thus the variance is given by
\begin{align*}
\variance{\spanlenrv{m}}
&= \expect{{\spanlenrv{m}}^2}-\expectsqr{\spanlenrv{m}}=\\
&= \frac{m\cdot(n+1)\cdot(mn+n+m)}{(m+2)(m+1)} - \left(\frac{m\cdot(n+1)}{m+1}\right)^2 = \\
&= \frac{m\cdot(n+1)}{m+1}\cdot\left(\frac{mn+n+m}{m+2} - \frac{m\cdot(n+1)}{m+1}\right) = \\
&= \frac{m\cdot(n+1)}{m+1}\cdot\left(\frac{mn+n+m}{m+2} - \frac{mn+m}{m+1}\right) = \\
&= \frac{m\cdot(n+1)}{m+1}\cdot\frac{(m+1)(mn+n+m)-(m+2)(mn+m)}{(m+2)(m+1)} =\\
&= \frac{m\cdot(n+1)}{m+1}\cdot\frac{n-m}{(m+2)(m+1)} =\\
&= \frac{m\cdot(n+1)\cdot(n-m)}{(m+1)^2\cdot(m+2)}.\qedhere
\end{align*}
\end{proof}

\subsection*{Acknowledgments}
The author is supported by an Adams Fellowship of the Israeli Academy of Sciences and Humanities.
This work was supported in part by ISF grant 230/10, by the Google Inter-University
Center for Electronic Markets and Auctions, and by
the European Research Council under the European Community's Seventh Framework
Programme (FP7/2007-2013) / ERC grant agreement no.\ [249159].
The author would like to thank his Ph.D.\ advisor, Sergiu Hart,
for useful discussions and comments,
and in particular for suggesting the idea underlying the third proof of \cref{distrib}(\labelcref{prob}).

\bibliographystyle{abbrv}
\bibliography{perms}

\end{document}